\documentclass{ws-procs9x6}

\def\End{{{\rm{End}}\,}}\def\hq{${\sf{HQ}}$}\def\id{{\rm{Id}}}\def\Aut{{\mbox{\rm Aut}\,}}\def\tr{{{\rm{tr}}}\,}
\def\R{\hat{R}}\def\F{\hat{F}}\def\rank{{{\rm{rank}}\,}}\def\trq{{{{\rm{tr}}}_{\scriptscriptstyle  \! \F}}}
\def\Ml{{{{\cal M}}_{\scriptscriptstyle SA}(\R,\F)}} 

\def\sa{{\scriptstyle SA}}

\def\for{\mathrm{for~}}\def\and{\mathrm{and~}}

\def\a{\alpha}\def\s{\sigma}\def\t{\tau}

\newcommand{\eod}{\hfill \rule{2.5mm}{2.5mm}}

\def\ok{\overline{k}}

\def\ba{\begin{eqnarray}}\def\ea{\end{eqnarray}}
\def\be{\begin{equation}}\def\ee{\end{equation}}\def\lb{\label}
\def\ni{\noindent}

\begin{document}
$\ $
\vskip -2cm

\title{Half-quantum linear algebra}

\author{A. Isaev}

\address{Bogoliubov Laboratory of Theoretical Physics, JINR\\
141980 Dubna, Moscow Region, Russia\\
E-mail: isaevap@theor.jinr.ru}

\author{O. Ogievetsky\footnote{On leave of absence from P. N. Lebedev Phys. Inst.,
Leninsky Pr. 53, Moscow, Russia}}

\address{Center of Theoretical Physics\\
Aix Marseille Universit\'e, CNRS, UMR 7332, 13288 Marseille, France\\
Universit\'e de Toulon, CNRS, UMR 7332, 83957 La Garde, France\\
E-mail: oleg@cpt.univ-mrs.fr}

\begin{abstract}
The Cayley--Hamilton--Newton theorem for half-quantum matrices is proven.
\end{abstract}

\keywords{\hspace{-.05cm}Yang--Baxter,\hspace{-.05cm} quantum matrix,\hspace{-.05cm} Cayley--Hamilton,\hspace{-.05cm} Newton identities}

\bodymatter

\vskip 1cm\ni
{\bf 1. Introduction}

\vskip .2cm\ni
The Cayley--Hamilton theorem in the linear algebra and the Newton identities in the theory of symmetric functions have a common predecessor,
discovered in \cite{IOP} and called there Cayley--Hamilton--Newton (CHN) identities. These are matrix identities, which naturally include the Cayley--Hamilton
identity and whose trace reproduces the Newton identities. The CHN identities admit a $q$-deformation, see Ref. \refcite{IOP}.  
Later\cite{IOP1}, the CHN identities were established for a larger class of quantum matrix algebras which includes both the algebra of functions on a quantum group 
of $GL$ type and the deformed universal enveloping algebra of $\mathfrak{gl}$ type. This class of algebras, defined with the help of a compatible pair $(\R,\F)$ of R-matrices
(see the precise definition in Section 2), can be described differently, as braided Hopf algebras universally coacting (in the braided sense) on two, left and right,  
quantum spaces. When the R-matrix $\F$ is the flip, this is the usual universally coacting Hopf algebra. 

There is, nowadays, a certain interest in the algebras universally coacting on just one quantum space, see \cite{CFR} and the literature cited there. We give the definition of 
the braided Hopf algebra $\Ml$ universally coacting (again in the braided sense) on one quantum space. The definition appeals again to a compatible pair $(\R,\F)$ of 
R-matrices. The matrix whose entries generate $\Ml$ is called {\it half-quantum} (\hq-)matrix. An example of \hq-matrices arises naturally in the theory of quantum groups
of $GL$ type with spectral parameters, see Section 2. 

The main result of the present article is the CHN theorem for the \hq-matrices 
defined by a pair $(\R,\F)$ where $\R$ is of Hecke type, that is, the spectral decomposition of $\R$ has two projectors, $S$ and $A$. The algebras $\Ml$ 
coact on the ``right" quantum space associated to $S$. The formal exchange of projectors leads to 
the statements about the algebras coacting on the ``left" quantum space, associated to $S$, so it is enough to investigate the algebras $\Ml$. 

The coefficients of the CHN identities are analogues of the sets of elementary and complete symmetric functions. Contrary to the algebra generated by the full quantum matrix, none of these sets is commutative for \hq-matrices. 
As for full quantum matrices, the CHN theorem for \hq-matrices uses modified matrix powers. 
We define four, in general, different powers, $M^{\rightarrow{k}}$ and $M^{\leftarrow{k}}$, 
of an \hq-matrix $M$; if $M$ is a full quantum matrix then $M^{\rightarrow{k}}=M^{\leftarrow{k}}$. 
As a direct consequence of the CHN theorem we obtain the Cayley--Hamilton theorem and Newton identities for \hq-matrices.  

{}The \hq-matrices can be defined in a more general setting, when one knows which eigen-projectors of an R-matrix belong to the symmetric and anti-symmetric part (see Section 3 in\cite{O} for a general discussion).
This is the case, in particular, of the R-matrices for simple Lie groups (in any representation). The CHN theorem for full quantum matrices is known\cite{OP} in the simplest situation, for the orthogonal and symplectic groups in the defining representation and it will be interesting to generalize it for the \hq-matrices.

\ni{\bf Notation.} Let $V$ be a finite-dimensional $\mathbb C$-vector space and $V^{\otimes M}:=
V\otimes V\otimes\dots$ ($M$ times, $M=0,1,2,\dots$; by convention, $V^{\otimes 0}:={\mathbb C}$). Let
$\mathfrak{A}$ be an associative algebra and $X$ an $\mathfrak{A}$-valued operator on $V$ (that is, the matrix 
elements of $X$ belong to $\mathfrak{A}$). We denote by $X_j$ the operator which
acts as $X$ on the $j$-th copy of the space $V$ in $V\otimes V\otimes\dots$ (and as the identity on other copies); for
an $\mathfrak{A}$-valued operator $Y$ on $V\otimes V$,  we denote by $Y_{j,k}$ the operator which acts as $Y$ on the $j$-th and $k$-th
copies of the space $V$ in $V\otimes V\otimes\dots$ (and as the identity on other copies) {\it etc}. 
We write sometimes $Y_{j} := Y_{j,j+1}$  for brevity. The identity operator on $V^{\otimes M}$ we
denote by $\id$ when the value of $M$ is clear from the context. We denote by $P\in \Aut(V\otimes V)$ the permutation
operator, $P(v\otimes w):=w\otimes v$. The operation of taking traces in copies of $V$ with the numbers $i_1\dots i_k$
is denoted by $\tr_{(i_1\dots i_k)}$. For an $\mathfrak{A}$-valued operator $Z$ on $V^{\otimes m}$ define the operator
$Z^{\uparrow k}$, $k=0,1,2,\dots$, in $V^{\otimes m+k}$ by $Z^{\uparrow k}:=\id^{\otimes k}\otimes Z$.
We use symmetric $q$-numbers, $j_q := q^{j-1}+q^{j-3}+\dots+q^{-j+1}$ for $j=0,1,\dots$ (by convention, the empty sum is 0).
Given a numerical $q$, we shall say that a positive integer $j$ is $q$-admissible if $k_q\neq 0 \ \ \text{for}\ k=1,2,\dots,j$.

\vskip .6cm\ni
{\bf 2. Generalized half-quantum matrices}

\vskip .2cm
\ni{\bf 2.1 Operators $\R$ and $\F$.} 
Generalized half-quantum matrix algebras are defined with the help of 
two operators $\R,  \F\in \Aut(V\otimes V)$ which satisfy, depending on a property in question, to a combination of conditions {\bf (i)}-{\bf (iv)} below.

\vskip .2cm
\ni {\bf (i)} $(\R,  \F)$ is a  {\em compatible pair} of R-matrices: 
\ba\lb{ybe}&&
\R_{12}\R_{23}\R_{12}=\R_{23}\R_{12}\R_{23}\,  ,\qquad \F_{12}\F_{23}\F_{12}=\F_{23}\F_{12}\F_{23}\, ,\\ 
\lb{compat}&&\R_{12}\F_{23}\F_{12}=\F_{23}\F_{12}\R_{23}\, ,\qquad \F_{12}\F_{23} \R_{12}=\R_{23}\F_{12}\F_{23}\, .\ea

\ni {\bf (ii)} The R-matrix $\R$ is of  {\it Hecke type}: its projector decomposition is 
$\R =qS^{(2)}-q^{-1}A^{(2)}$ where $q\in\mathbb{C}^*$ is such that 2 is $q$-admissible. The projectors $S^{(2)}$ and $A^{(2)}$ 
(called, respectively, $q$-symmetrizer and $q$-antisymmetrizer)
are complementary, $S^{(2)}+A^{(2)}=\id$.

The higher $q$-symmetrizers and $q$-antisymmetrizers $S^{(k+1)},A^{(k+1)}\in \End(V^{\otimes k+1})$, $k>1$, can be defined, if $k+1$ is $q$-admissible, 
in several ways, {\it e.g.}, inductively, in any of the forms
\ba\lb{hsas1}  &&\,A^{(k)\uparrow 1}\left(q^{k}\,\id-k_q\R_{1}\right)A^{(k)\uparrow 1}=:(k+1)_qA^{(k+1)}:=\,A^{(k)}\left(q^{k}\,\id-k_q\R_{k}\right)A^{(k)}\ ,\\
\lb{hsas2} &&\,S^{(k)\uparrow 1}\left(q^{-k}\,\id+k_q\R_{1}\right)S^{(k)\uparrow 1}=:(k+1)_qS^{(k+1)}:=\,S^{(k)}\left(q^{-k}\,\id+k_q\R_{k}\right)S^{(k)}
\ .\ea
By convention, $A^{(1)}:=\id$ and $S^{(1)}:=\id$; then (\ref{hsas1}) and (\ref{hsas2}) hold for $k=1$ as well.

\ni {\bf (iii)} The Hecke R-matrix $\R$ is {\it even of height\cite{G} $n$}: $n$ is $q$-admissible,
$\rank (A^{(n)}) = 1$ and $A^{(n)}(q^{n}\,\id-n_q\R_{n})A^{(n)}=0$.

\ni {\bf (iv)} The R-matrix $\F$ is invertible and {\it skew-invertible}: there exists
$\Psi\in \Aut(V\otimes V)$ satisfying $\tr_2(\Psi_{12}\F_{23})=P_{13}$.

\ni{\bf 2.2 F-trace.} The {\em quantum trace} of a matrix $X$ with arbitrary entries is 
$\trq (X):=\tr (D X)$, where $D\in \End(V)$ is defined by $D_1:=\tr_2(\Psi_{12})$. We have
\[ \trq_{(2)} (\F_1) = \id_1\ ,\ \F_1 D_1 D_2 = D_1 D_2 \F_1\ ,\ 
\trq_{(2)} (\F_1^{\epsilon} \, X_1 \, {\F_1}^{-\epsilon}) = \id_1 \, \trq (X)\ ,\ \epsilon=\pm1\ .\]

\ni{\bf 2.3} Let $M$ be an operator on $V$ with arbitrary entries. Define inductively 
\be\lb{Mk}M_{\overline{1}} := M_1\ , \qquad M_{\overline{k+1}}:= \F_{k, k+1}^{\phantom{-1}} M_{\ok} \F_{k, k+1}^{-1}\ .\ee

\ni{\bf Definition.} {\it Let $\Ml$ be the unital associative algebra generated by the entries of $M$ with the defining relations
\be\lb{QMAl}S^{(2)}_{12}\, M_{\overline{1}} \, M_{\overline{2}}\,  A^{(2)}_{12} = 0 \ .\ee

\noindent We call $M$ the half-quantum $\sa$-matrix (\hq-matrix, for brevity).}

To explain the meaning of the definition, let ${\cal{X}}_{S}$ be the algebra generated by components $\{ x^a\}$ of a vector from $V$ with the defining quadratic relations 
$(S^{(2)})^{ab}_{cd}x^cx^d=0$ (the Einstein summation convention is assumed throughout the text). The \hq-matrix coacts, $x^a\to M^a_b x^b$, on the algebra ${\cal{X}}_{S}$ and generates the  
braided Hopf algebra universally coacting on ${\cal{X}}_{S}$. The braiding between the entries of $M$ and $x$ is given by $X_1M_2=
M_{\overline{2}}X_1$, where $X^a_b:=c_bx^a$, $c_b$ are independent commuting variables. 
As for the braided coproduct, let $M$ and $M'$ be \hq-matrices which satisfy
$$\F_{12}^{\phantom{1}}M_1^{\phantom{1}} \F_{12}^{-1}{M'}_1^{\phantom{1}}=
{M'}_1^{\phantom{1}} \F_{12}^{\phantom{1}}M_1^{\phantom{1}} \F_{12}^{-1}\ .$$
Then $MM'$ is again an \hq-matrix.

\ni{\bf Notation.} 
{}For natural $l$ and $k$, $l\leqslant k$, set $M_{\overline{l\to k}}:=M_{\overline{l}} \dots M_{\ok}$, $\F_{l\to k}:=\F_{l} \dots \F_{k-1}$, 
$\R_{l\to k}:=\R_{l} \dots \R_{k-1}$ and $\R_{k\leftarrow l}:=\R_{k-1}\dots \R_{l}$ (by convention, the empty product is 1).

\ni{\bf 2.4}  Let $(\R,  \F )$ be a pair satisfying the conditions {\bf (i)} and {\bf (iv)}, Subsection {\bf 2.1}. The following Lemma, valid for
a matrix $M$ with arbitrary entries, was proved in Ref. \refcite{IOP1}. 

\ni {\bf Lemma 1.~}{\em {\bf a)~} We have
\ba\lb{l1}&&\F_i M_{\ok} = M_{\ok} \F_i\ \ \and \ \R_i M_{\ok}= M_{\ok} \R_i\ \ \for  k\neq i, i+1\ ,\\
\lb{l2}&&\F_{i\rightarrow k+1}M_{\overline{i\to k}}=M_{\overline{i+1\to k+1}}\F_{i\rightarrow k+1}\ \ \for  i\leq k\ .\ea
\ni{\bf b)~} Let $\a(Y^{(k)}):=\trq_{(1,\dots ,k)}(Y^{(k)}M_{\overline{1\to k}})$ where $Y^{(k)}$ is a polynomial in $\R_1,\dots,\R_{k-1}$. Then
\be\lb{l4}\trq_{(i+1,\dots ,i+k)}(Y^{(k)\uparrow i} M_{\overline{i+1\to i+k}}) = \id_{1,\dots ,i}\, \a(Y^{(k)})\  .\ee}

\ni{\bf 2.5} Assume now that $M$ is an \hq-matrix. Since $\R_1=qS^{(2)} -q^{-1}A^{(2)}$ and $S^{(2)} + A^{(2)} = \id$,
the relations (\ref{QMAl})  can be rewritten in the four following equivalent
(when $q+q^{-1}\neq 0$) forms
\be \lb{hq02a}\begin{array}{l}\R_1 M_{\overline{1}} M_{\overline{2}} A^{(2)} = - q^{-1} \, M_{\overline{1}} M_{\overline{2}} A^{(2)} \; ,
\;\;\;  S^{(2)} M_{\overline{1}} M_{\overline{2}} \R_1 = q \, S^{(2)} M_{\overline{1}} M_{\overline{2}}  \; ,\\[.4em]
A^{(2)} M_{\overline{1}} M_{\overline{2}} A^{(2)} =M_{\overline{1}} M_{\overline{2}} A^{(2)} \; , \;\;\;
S^{(2)} M_{\overline{1}} M_{\overline{2}} S^{(2)}=  S^{(2)} M_{\overline{1}} M_{\overline{2}}  \; . \end{array}\ee

\ni By induction on $k=1,2,\dots$ we obtain, using Lemma 1,
\ba\lb{l3} &S^{(2)}_{j,j+1} M_{\overline{j}}^{\phantom{1}} M_{\overline{j+1}}^{\phantom{1}} A^{(2)}_{j,j+1} = 0 \ ,&\\[.2em]
\lb{hq03}&A^{(k)\uparrow i} M_{\overline{i+1\to i+k}} A^{(k)\uparrow i} = M_{\overline{i+1\to i+k}} A^{(k)\uparrow i}\ ,\ i=0,1,2,\dots& \\[.4em]
\lb{hq03b}&S^{(k)\uparrow i} M_{\overline{i+1\to i+k}} S^{(k)\uparrow i}=S^{(k)\uparrow i}  M_{\overline{i+1\to i+k}}\ ,\ i=0,1,2,\dots&\ea
{\ni}In eqs.(\ref{hq03})-(\ref{hq03b}), $k$ is assumed to be $q$-admissible. 

\ni{\bf 2.6 Example.} Let $(\R,\F)$ be a compatible pair with $\R$ of Hecke type. Let $\R(u):=\R+(q-q^{-1})/(u^2-1)$ be the Baxterizaton of $\R$ and let $T(u)$ be the matrix 
generating the algebra with the defining relations $\R_{12}(u/v)T_{\overline{1}}(u)T_{\overline{2}}(v)=T_{\overline{1}}(v)T_{\overline{2}}(u)\R_{12}(u/v)$ (in the simplest 
cases, this can be an algebra related to quantum or classical (super)-Yangians or affine algebras). Now, $\R(q)$ is proportional to the symmetrizer, so setting $u=qv$, multiplying from the right by $A^{(2)}_{12}$ and writing $T(qv)=q^{v\partial_v}T(v)q^{-v\partial_v}$ we find
$S^{(2)}_{12}q^{v\partial_v}T_{\overline{1}}(v)q^{-v\partial_v}T_{\overline{2}}(v)A^{(2)}_{12}=0$,
which implies that $q^{-v\partial_v}T$ and $Tq^{-v\partial_v}$ are \hq-matrices, $S^{(2)}_{12}q^{-v\partial_v}T_{\overline{1}}(v)q^{-v\partial_v}T_{\overline{2}}(v)A^{(2)}_{12}=0$
and $S^{(2)}_{12}T_{\overline{1}}(v)q^{-v\partial_v}T_{\overline{2}}(v)q^{-v\partial_v}A^{(2)}_{12}=0$.

\vskip .6cm\ni
{\bf 3. Cayley--Hamilton--Newton theorem}

\vskip .2cm
\ni{\bf 3.1 Symmetric functions.} Let $M$ be an \hq-matrix.   Define four sets, $\{ s_k(M) \}$,
$\{ \overline{s}_k(M)\}$, $\{ \s_k(M)\}$ and $\{ \tau_k(M) \}$, $k=0,1,2,\dots$, of elements of the algebra $\Ml$ 
by $s_0(M) =\overline{s}_0(M) =\s_0(M) = \t_0(M) = 1$ and, for $k>0$, 
\ba\lb{s}&
s_k(M) :=\trq_{(1\dots k)}(\R_{1\to k} M_{\overline{1\to k}})\ \ ,&\ \ \overline{s}_k(M) :=\trq_{(1\dots k)}(\R_{k\leftarrow 1} M_{\overline{1\to k}})\, ,\\[.2em]
\lb{sig} &
\s_k(M) :=q^k \, \trq_{(1\dots k)}(A^{(k)}M_{\overline{1\to k}})\ \ ,&\ \  
\t_k(M) := q^{-k} \, \trq_{(1\dots k)}(S^{(k)}M_{\overline{1\to k}})\, .\ea

\ni{\bf 3.2 Powers of \hq-matrices.}
Define the matrices $M^{\rightarrow{k}}$, $\tilde{M}^{\rightarrow{k}}$, $M^{\leftarrow{k}}$ and $\tilde{M}^{\leftarrow{k}}$, {\em $k$-th powers} of $M$, by
\ba\lb{power}&(M^{\rightarrow{k}})_1:= \trq_{(2,\dots ,k)}(M_{\overline{1\to k}}\R_{1\to k})\ ,\  
(M^{\leftarrow{k}})_1:= \trq_{(2,\dots ,k)}(\R_{k\leftarrow 1}M_{\overline{1\to k}})\ ,&\\
\lb{power2}&(\tilde{M}^{\rightarrow{k}})_1:= \trq_{(2,\dots ,k)}(\R_{1\to k}M_{\overline{1\to k}})\ ,\  
(\tilde{M}^{\leftarrow{k}})_1:= \trq_{(2,\dots ,k)}(M_{\overline{1\to k}}\R_{k\leftarrow 1})\ .&\ea
Define the $k$-th wedge powers $M^{{\wedge k}}$, $\tilde{M}^{{\wedge k}}$ and the 
$k$-th symmetric powers $M^{{{\scriptscriptstyle\cal S} k}}$, $\tilde{M}^{{{\scriptscriptstyle\cal S} k}}$ of $M$ by 
\ba\lb{skew-power}&(M^{{\wedge k}})_1:= \trq_{(2,\dots ,k)}(A^{(k)} M_{\overline{1\to k}})\ ,\  
(M^{{{\scriptscriptstyle\cal S} k}})_1:= \trq_{(2,\dots ,k)} ( M_{\overline{1\to k}} S^{(k)})\ ,&\\
\lb{skew-power2}&(\tilde{M}^{{\wedge k}})_1:= \trq_{(2,\dots ,k)}(M_{\overline{1\to k}}A^{(k)} )\ ,\  
(\tilde{M}^{{{\scriptscriptstyle\cal S} k}})_1:= \trq_{(2,\dots ,k)} (S^{(k)} M_{\overline{1\to k}})\ .&\ea

\ni{\bf 3.3 Cayley--Hamilton--Newton theorem for \hq-matrices} (Cf. \cite{IOP}) {\em Assume that the pair $(\R,  \F)$ satisfies conditions 
{\em {\bf (i)}, {\bf (ii)}} and {\em {\bf (iv)}}, Subsection {\bf 2.1}, and assume that the positive integer $j$ is $q$-admissible. The \hq-matrix satisfies the following matrix identities}
\ba\lb{chn}&\displaystyle{j_q\, M^{{\wedge j}}=\sum_{k=0}^{j-1}(-1)^{j-k+1}M^{\leftarrow{j-k}}\s_k(M)}\ ,&\\
&\displaystyle{j_q M^{{{\scriptscriptstyle\cal S}j}}=\sum_{k=0}^{j-1} M^{\rightarrow{j-k}}\ \t_k(M)\, ,}&\\
\lb{chn2}&\displaystyle{j_q\, \tilde{M}^{{\wedge j}}=\sum_{k=0}^{j-1}(-1)^{j-k+1}\tilde{M}^{\leftarrow{j-k}}\s_k(M)}\ ,&\\ 
&\displaystyle{j_q \tilde{M}^{{{\scriptscriptstyle\cal S}j}}=\sum_{k=0}^{j-1} \tilde{M}^{\rightarrow{j-k}}\ \t_k(M)\, .}&\ea

\ni We sketch the {\it proof} of (\ref{chn}). For $k$, $1\leq k<j$, rewrite $(M^{\leftarrow{j-k}}\s_k(M) )_1$ as
\be\lb{dlass}\begin{array}{c} q^{k}\trq_{(2\dots j-k)}(\R_{j-k\leftarrow 1} M_{\overline{1\to j-k}})
\trq_{(1\dots k)}(\! A^{(k)}\! M_{\overline{1\to k}})
\\[1em]
=q^k\trq_{(2\dots j)}(\! A^{(k)\uparrow j-k}\! \R_{j-k\leftarrow 1}M_{\overline{1\to j-k}} M_{\overline{j-k+1\to j}})
=q^k\trq_{(2\dots j)}(\! A^{(k)\uparrow j-k}\! \R_{j-k\leftarrow 1}M_{\overline{1\to j}})\ .\end{array}\ee
{\ni} Here we used (\ref{l4}). We use the recurrence for the $q$-antisymmetrizers in the form
$q^kA^{(k)\uparrow 1}=(k+1)_qA^{(k+1)}+k_qA^{(k)\uparrow 1}\R_1 A^{(k)\uparrow 1}$
to rewrite the last expression in (\ref{dlass}) as
\be\lb{dlass2}\begin{array}{l}(k+1)_q\trq_{(2\dots j)}(A^{(k+1)\uparrow j-k-1}\R_{j-k\leftarrow 1}M_{\overline{1\to j}})\\[.6em]
\hspace{1cm}+k_q\trq_{(2\dots j)}(A^{(k)\uparrow j-k}\R_{j-k}A^{(k)\uparrow j-k}\R_{j-k\leftarrow 1} M_{\overline{1\to j}})\ .\end{array}\ee
\vskip -.1cm{}
{\ni}In the second term of (\ref{dlass2}), the right $A^{(k)\uparrow j-k}$ commutes with $\R_{j-k\leftarrow 1}$ and the left $A^{(k)\uparrow j-k}$ can be cyclically moved. So, the
last term takes the form
\[\begin{array}{l} k_q\trq_{(2\dots j)}(\R_{j-k+1\leftarrow 1} A^{(k)\uparrow j-k} M_{\overline{1\to j}} A^{(k)\uparrow j-k})
= k_q\trq_{(2\dots j)}(\R_{j-k+1\leftarrow 1} M_{\overline{1\to j}} A^{(k)\uparrow j-k})\\[.6em]
\hspace{1cm}= k_q\trq_{(2\dots j)}(A^{(k)\uparrow j-k}\R_{j-k+1\leftarrow 1}M_{\overline{1\to j}})\ .\end{array}\]
\vskip -.1cm{}
{\ni}We used (\ref{hq03}) in the first equality; in the second equality, we moved
$A^{(k)}$ to the left using the cyclic property of the trace. Finally, we obtain
\be\lb{hq11}\begin{array}{rcl}
(M^{\leftarrow{j-k}}\s_k(M) )_1&=&(k+1)_q\trq_{(2\dots j)}(A^{(k+1)\uparrow j-k-1}\R_{j-k\leftarrow 1}M_{\overline{1\to j}})\\[.6em]
&&+ k_q\trq_{(2\dots j)}(A^{(k)\uparrow j-k}\R_{j-k+1\leftarrow 1}M_{\overline{1\to j}})\ .\end{array}\ee
{\ni}Eq.(\ref{chn}) follows from the expressions (\ref{hq11}) for $M^{\leftarrow{j-k}}\s_k(M)$, $k=1,2,\dots ,j-1$, and 
$M^{\leftarrow j}\s_0(M)=M^{\leftarrow j}$.\eod

\ni{\bf 3.4 Corollaries.}

\ni{\bf Newton relations.} Taking the trace in (\ref{chn}), we find, for a $q$-admissible $j$, 
\ba\lb{nr1} q^{-j}j_q \s_j(M)=\sum_{k=0}^{j-1}(-1)^{j-k+1}\  \overline{s}_{j-k}(M) \s_k(M)\ ,\\
\lb{nr2}q^j j_q \t_j(M)=\sum_{k=0}^{j-1} s_{j-k}(M) \t_k(M)\ .\ea

\ni{\bf Cayley-Hamilton theorem.} Assume that $\R$ is even of height $n$. Then there exist tensors $\epsilon_{a_1a_2\dots a_n}$ and $\epsilon^{a_1a_2\dots a_n}$ 
such that 
$$(A^{(n)})_{a_1\dots a_n}^{b_1\dots b_n}=\epsilon_{a_1\dots a_n}\epsilon^{b_1\dots b_n}\ \ \text{and}\ \ \epsilon_{a_1\dots a_n}\epsilon^{a_1\dots a_n}=1\ .$$ 
By (\ref{hq03}) for $k=n$ and $i=0$, we have 
\be\lb{deqde}(M_{\overline{1\to n}})^{a_1\dots a_n}_{b_1\dots b_n}\epsilon^{b_1\dots b_n}={\text{det}}_q(M)\epsilon^{a_1\dots a_n}\ ,\ee
where the {\it quantum determinant} of the \hq-matrix $M$ is defined, up to a normalization factor, by 
$${\text{det}}_q(M):=\epsilon_{a_1\dots a_n}(M_{\overline{1\to n}})^{a_1\dots a_n}_{b_1\dots b_n}\,\epsilon^{b_1\dots b_n}\ .$$ 
Eq.(\ref{deqde}) implies that 
$$ \tilde{M}^{{\wedge n}}={\text{det}}_q(M) {\cal{D}}\ ,$$
where the matrix ${\cal{D}}$ is defined by 
$$ {\cal{D}}_1:=\trq_{(2\dots n)}(A^{(n)})\ .$$
Taking $j=n$ in (\ref{chn2}) we obtain the Cayley--Hamilton theorem for the half-quantum matrices:
\be\label{CHNth} \sum_{k=0}^{n-1}(-1)^{n-k}\tilde{M}^{\leftarrow{n-k}}\s_k(M)+
n_q{\text{det}}_q(M){\cal{D}}=0\ .\ee

\ni{\bf Standard $\R$.} Let $(\R_{DJ},P)$ be the compatible pair with the standard multi-parametric Drinfeld--Jimbo R-matrix in dimension $n$. 
In \cite{Z}\,, a version of the Cayley--Hamilton theorem, with usual matrix powers but diagonal matrices as coefficients, has been proven. Namely, the \hq-matrix $M$ 
satisfies the matrix identity 
\begin{equation}\label{diagcoeff} M^n+\sum_{k=1}^n (-1)^k\Sigma_kM^{n-k}=0\, ,\end{equation}
{\ni}with certain diagonal matrices $\Sigma_m$, $m=1,\dots,n$.  
By the same technics as in \cite{OV} one can relate modified
powers of $M$ - in this situation, i.e., for the compatible pair $(\R_{DJ},P)$ - with usual powers and deduce 
(\ref{diagcoeff}) from (\ref{CHNth}).

\end{document}